\documentclass[12pt]{article}
\usepackage[dvips]{graphics}
\usepackage[pdftex]{graphicx}
\usepackage{graphicx}

\usepackage{amsmath}
\usepackage{amssymb}
\usepackage{amsthm}
\usepackage{epsfig}
\usepackage{amscd}
\usepackage{setspace}
\usepackage{array}
\usepackage{multirow}
\usepackage{anysize}
\usepackage{latexsym}
\usepackage{verbatim}

\usepackage{color}
\usepackage{epsfig}
\usepackage{amssymb}

\label{marker}

\title{\bf{ Poincar\'{e}  map construction for some classical two predators - one prey systems.}}

  \author{A.\,V.\,Osipov \\
  Saint-Petersburg University, Russia \\
  email: av\_osipov@mail.ru  \\
    G.\,J.\,S\"oderbacka  \\
    \AA bo Akademi, Finland \\
    email: gsoderba@abo.fi}

\date{}

\begin{document}
\maketitle

\begin{abstract}

In this paper we consider a family of system with 2 predators feeding on one prey. We show how to construct a positively invariant set in which it is possible to define a Poincar\'{e} map for examining the behaviour of the system, mainly in the case when both predators survive. We relate it to examples from earlier works.

\end{abstract}

\section{Introduction}

 The question of coexistence in models of several predators feeding on one prey have been considered in many works. One question is to find parameter criteria for coexistence. A complicated question is to find out the form of the coexistence. Mathematically it means to examine the structure of the attractor, which can be periodic or of different chaotic types. In this work we introduce a construction of a positively invariant set, where it is possible to examine the attractor by estimating properties for a  Poincar\'{e} map. We also find out when it is possible or not to find the invariant set.

 We consider the system
 \begin{equation}\label{a1}
 \dot{x_1}={\phi}_1(s) x_1, \quad  \dot{x_2} = {\phi}_2(s)x_2,  \quad\quad
 \dot{s}=H(x_1,x_2,s).
 \end{equation}
 The function  $H = H(x_1,x_2,s)$  has the form
 $H=h(s)-{\psi}_1(s)x_1-{\psi}_2(s)x_2 .$ \\
 \vskip0mm
We assume the following conditions:
 \vspace{1mm}

 \noindent
 $A_1 :\quad$ All functions belong to the class 
 $\, C^2 [ 0,\infty ) ,\,$
 and the variables $x_1, x_2, s$ are in the octant 
  $\, R_3^+ \colon \;
 x_1 \geq 0 , x_2 \geq 0, s \geq 0 ,\,$ which is an invariant set.
 \vskip0.5mm

 \noindent
 $A_2 : \quad   {\psi}_i(0)=0, \quad {\psi}_i'(s)>0\,
 $ for  $s>0$.
 \vskip0.5mm

 \noindent
 Here and further we will suppose, that  $i$  takes values from the set $\{1,\,2\}. $
 \vskip0.5mm

\noindent $A_3 : \quad   {\phi}_i'(s)>0\,
$ for  $s>0,\quad$ and there exists 
 ${\lambda}_i > 0$ such that  ${\phi}_i({\lambda}_i)=0. \quad$
\vskip0.5mm

\noindent $A_4 : \quad  h(0) = h(1) = 0 , \quad  h'(0) > 0 \;$  and
 $\, h''(s) < 0 \;$ for $ s>0. $
\vskip0.5mm

\noindent $A_5 : \quad   {\psi}_1 (s) < {\psi}_2 (s) \, $ for
$s>0$. \vskip0.5mm

\noindent $A_6 : \quad $ There exists  $\tau > 0\,$ such that
 $\, {\phi}_1(\tau) = {\phi}_2(\tau)  $
and
  ${\phi}_2(s)> \phi_1 (s) \; \Leftrightarrow \; s > \tau . $
\vskip0.5mm

 \noindent
 $A_7 : \quad  0  <  \tau  <  {\lambda}_2 <  {\lambda}_1  <  1 .  $
 \\

\vskip0mm

 The vector field of system (\ref{a1}) we will denote by $\,{\cal F}. $
\vskip1mm

    Systems with the following special choice of
  functions

 \begin{equation}\label{a2}
 h(s)=\gamma s(K-s),\qquad {\psi}_i(s)={\alpha}_i\frac{s}{s+a_i},\qquad
 {\phi}_i(s)=m_i{\psi}_i(s)-d_i,
 \end{equation}
 where all introduced parameters are positive, have been examined extensively.
 The general behaviour of such systems has been studied in 
 \cite{Hsu4,Hsu3} and parameter conditions for different behaviour are given in 
 \cite{Hsu5}. The possibility for coexistence connected with bifurcations were examined in \cite{smith1,ButlerWalt,Keener}. Many works have considered the situation when function $h$ is linear and thus condition $A_4$ is violated.
 More general systems are considered in \cite{smiththieme,alebra}.

 Such systems can have chaotic regimes. In particular, in 
 \cite{osipov3,osipov4,osipov5} was given a bifurcation diagram with a chain of period doubling bifurcations for 
  system (\ref{a1}) with (\ref{a2}), where 
 \begin{equation}\label{a3}
 h(s)=s(1-s),\qquad {\phi}_i(s)=\frac{s-{\lambda}_i}{s+a_i},\qquad
 {\psi}_i(s)=\frac{s}{s+a_i}.
 \end{equation}
 The system with functions (\ref{a2}) were transformed into this system
 in \cite{Keener}. In addition in \cite{osipov3,osipov4,osipov5}  was assumed $m_i=1$ for these new parameters after the transform. 
  A  one-dimensional discrete model was built with the same chain
 of bifurcations as  for the real Poincar\'{e} map.
 It was also shown that along this bifurcation path there is the possibility that  
  the Poincar\'{e} map
 cannot be properly defined and there can appear a so called "spiral-like chaos".
 Here we modify the bifurcation path to correspond to our technics of constructing the invariant set.
 \smallskip

 In this work we consider a broad class of systems of the type
 (\ref{a1}),
 including the system with the right hand sides of form
 (\ref{a2}) and with  (\ref{a3})
 as the main example .
 We formulate sufficient conditions for correct definition of 
 Poincar\'{e} map and, as a corollary, for absence of "spiral-like chaos".
 
 We note that it is well known that such systems are not directly applicable in biology without modifications, for example, because in the model some populations can go unacceptable low and still survive. Also in recent time most models include more complicated dependences on the variables. Usually the equations for predators are not fully independent on the other predator and functions $\psi$ and $\phi$ can include such dependence. Anyhow all these modifications often essentially use the technics worked out for our systems, why the study of system (\ref{a1}) is still motivated.
 
The outline of the paper is the following. After this introduction we give a short overview of local behaviour of the system. In Section 2 we consider the two dimensional behaviour in the  $(x_i, s)$-coordinate planes. Section~3 is the main section where we introduce comparison systems and the concept of a tangency curve we use for constructing the positively invariant set.  The boundary of the set will mainly be formed by parts of trajectories of a glued system. Finally in Section 4 we give some numerical results for cases, when the construction is possible and when not possible containing  more complicated chaos.

\subsection{Equilibria}

 From conditions  $A_1-A_7$  
 follows, that the system has four equilibria:

\noindent
 $1^{\circ}. \quad {\bf O}\,(0,0,0), $  which is saddle with two-dimensional stable manifold $\,s=0\,$   and one-dimensional unstable: $ x_1 = x_2 = 0, \; s<1. $

\noindent $2^{\circ}. \quad {\bf O_1}\,(0,0,1), $   
 which is a saddle 
 with one-dimensional stable manifold
 $\, x_1 = x_2 =
0\,$
and a two dimensional unstable one.

 \noindent
 $3^{\circ}.\quad$ Two points in the coordinate planes:
 
 $\; {\bf P_1}\,(x_1^0, 0, \lambda_1)\,$ и
 $\; {\bf P_2}\,(0, x_2^0, \lambda_2), $
 with index 1 in these coordinate planes.

 \smallskip
 From conditions $A_3$ and $A_7$ follows that there are never any equilibria
 in the interior of the positive octant.

 The type of the equilibria are immediately deduced from the form
 of the Jacobian matrix.
 At the equilibrium  $ {\bf O^*}\,(x_1^*, x_2^*, s^*)\,$
 the Jacobian matrix has the form
 \begin{equation}
  J^* \; =\; J\,({\bf O^*} )\quad  = \quad
  \left(
 \begin{array}{ccc}
 \phi_1 ( s^* ) & 0 &  \phi_1' ( s^* ) x_1^* \\
           0  &  \phi_2 ( s^* ) &   \phi_2' ( s^* ) x_2^* \\
   -\psi_1 ( s^* )  &  - \psi_2 ( s^* ) &   h'(s^*) -
 \psi_1' ( s^* ) x_1^* - \psi_2' ( s^* ) x_2^*     
  \end{array} \right) \;.
 \end{equation}

Observing that the sign of the real parts of the eigenvalues are
determined by the sign of $\phi_i (s^ *)$ 
and $   h'(s^*) -
 \psi_1' ( s^* ) x_1^* - \psi_2' ( s^* ) x_2^*$
we can find the type of the equilibria of ${\bf P_1}$  and ${\bf P_2}$. 
${\bf P_1}$ is a saddle with two dimensional stable manifold 
in the coordinate plane, if
$ h'(\lambda_1) -\psi_1' ( \lambda_1 ) x_1^0 <0  $, and a source, if $ h'(\lambda_1) -
 \psi_1' ( \lambda_1 ) x_1^0 >0  $.
 ${\bf P_2}$ is stable, if
$ h'(\lambda_2) -\psi_2' ( \lambda_2 ) x_2^0 <0  $, and a saddle with one dimensional stable manifold, if $ h'(\lambda_2) -
 \psi_2' ( \lambda_2 ) x_2^0 >0  $.

  \section{Limit cycles in the coordinate planes}
\vskip1mm

The coordinate planes of system (1) are invariant two dimensional systems
and in this section we describe some known behaviour of these systems.

\noindent
Here we use a version of 
 Zhang Zhi Fen 
 theorem \cite{Zeng,ZhangZhiFen,ZhangZhiFen2}:
\vskip1mm

{\bf Theorem 1.}\quad {\em Consider the system
$$ \dot{x} = - \varphi ( y ) -  F ( x ) , \quad
 \dot{y} = g ( x )  . $$
Assume:

\noindent 1)\quad $g(x)\;$
satisfies the local Lipschitz condition,  
  $xg(x)>0,\;x\neq 0,\quad
 G(+\infty)=G(-\infty)=+\infty,$  где   $G(x)=\int^x_0 g(x) dx;$

\noindent 2)\quad  $F'(x)$ is continuous,  $F'(x)/g(x)$ is non-decreasing on intervals
      $( - \infty , 0 )$ and $( - \infty , 0 )$,
  and is non-constant in any neighbourhood of       
       $x=0. $  
       
\noindent 3)\quad $\varphi (y)\;$ 
satisfies the local Lipschitz condition,
  $ y \varphi (y) > 0,\; y \neq 0, \qquad  \varphi (y) \,$
 is non-decreasing, $\quad
 \varphi (-\infty)=-\infty , \quad \varphi (+\infty)= +\infty \,;
 \quad \varphi (y)\,$ 
  has left and right derivatives at the point 
 $x=0,$ 
  which are non-zero in the case 
 $F'(0)=0. $
 
 Then the system does not have more than on limit cycle
 and, if it exists, it is stable. 
 }\\

In  \cite{Sunhong} 
was shown, that the problem of limit cycles in two dimensional systems
of type "predator-prey" often is solved by applying the 
 Zhang Zhi Fen theorem.

In particular, this theorem is possible to use for the system

$$\dot{s} = h(s) -\psi (s) x,\quad \dot{x}=\phi(s) x. \eqno(*)$$

In this case the system is transformed to the form

$$ \dot{s} = F(s) -  x  , \quad
 \dot{x} = \varphi ( s )x  . $$

 where $F(s)=h(s)/\psi(s), \quad \varphi(s)=\phi(s)/\psi(s).$

We formulate the corresponding theorem in
 \cite{Sunhong}, 
 applied to the two dimensional restrictions of system (\ref{a1})
 to the coordinate planes.
 In the proof the system is transformed to a Lienard system used in Theorem 1.

{\bf Theorem 2.}\quad {\em Consider the system
$$ \dot{s} =  F ( s ) -x , \quad
 \dot{x} = \varphi ( s )x  , $$
 where $F(s)=h(s)/\psi(s), \quad \varphi(s)=\phi(s)/\psi(s).$

Suppose:

\noindent 1)\quad All functions are continuous and differentiable on
 $(0,\infty).$

\noindent 2)\quad  There exists  $\lambda  >0$ such that $\varphi(\lambda )=0$ and
$(s-\lambda )\varphi(\lambda ) >0, s\neq \lambda.$

\noindent 3)\quad $F'(s)/\varphi(s)$
is non-increasing on the intervals
  $(0,\lambda)$ and
$(\lambda,1)$.

Then the system has no more than one limit cycle which is always stable.
 }

\bigskip
We now consider the standard version of the two-dimensional system
corresponding to the functions in (\ref{a3}):
$$\dot{x}=\frac{s-\lambda}{s+a} x,\quad \dot{s}= s(1-s-\frac{x}{s+a}).$$


By scaling time it is transformed to the system
$$\dot{x}=(s-\lambda) x,\quad \dot{s}= s((1-s)(s+a)- x).$$

In the first quadrant the system has an equilibrium
$$P((1-\lambda)(\lambda +a),\, \lambda),$$
which is an asymptotically stable node or focus for
$\lambda \geq \frac{1-a}{2}$ 
and repelling node or focus for 
 $$\lambda <\frac{1-a}{2}.$$ 
Using Theorem 2 this condition implies there will be a unique limit cycle in the plane.

\smallskip

Our main interest is in the case, where there are two limit cycles in
both coordinate planes. This requires
$ h'(\lambda_i) -
 \psi_2' ( \lambda_i ) x_i^0 >0  $ in the general case
 and $\lambda_i <\frac{1-a_i}{2}$ in our special case.

\section{The Poincar\'{e} map}

While the behaviour of the two-dimensional systems in the coordinate
plane is well known, the behaviour outside the coordinate planes can be very 
complicated in the case both predators survive and we can find behaviour not yet well understood.

It is clear from assumptions that all interesting behaviour will be
in the region where $s<1$.

In the system when the functions are chosen as in (\ref{a2}) and (\ref{a3}), 
if $a_1>a_2$ in \cite{kustarov} is shown that the predator $x_1$ goes extinct if

$$\gamma > 1 +\frac{a_1}{\lambda_1} \left( \lambda_1 -\lambda_2\right)
+ \frac{\lambda_2}{a_1} \left( a_1 -a_2\right), \,\,
\gamma =\frac{\lambda_2 a_1}{\lambda_1 a_2} $$

\noindent or equivalently if

$$\lambda_1 > \frac{a_1\lambda_2 (a_2+1)}{a_1a_2 +\lambda_2 (a_1-a_2) +a_2}  $$

and the predator $x_2$ goes extinct in the case $\lambda_1<\lambda_2$.
Thus in this case we will concentrate on the case when these inequalities are not satisfied.

 \subsection{Decomposition. Simple change of variables}

In order to construct a positively invariant set and understand something
about the behaviour inside this set we use comparison systems with simpler
two-dimensional behaviour.

In order to find the comparison systems we consider a transform of
the original system into new variables:
  $  \  m= x_1 + x_2 ,
 \quad \xi = x_1 / m  . $
 In the new variables system (\ref{a1}) obtains the form:
\begin{equation}\label{5}
\dot{m}= p ( s,\,\xi ) \,m,   \quad\quad \dot{\xi} = \sigma ( s ) \,
\xi ( 1 - \xi ),   \quad\quad \dot{s}= h(s)- q (s, \,\xi )\,m,
\end{equation}
 where
\begin{equation}\label{b2}
 p ( s,\, \xi ) = {\phi}_1(s) \xi + {\phi}_2(s)( 1 - \xi ),   \quad\quad
 \sigma ( s ) = {\phi}_1(s)  - {\phi}_2(s) ,   \quad\quad
 q (s,\, \xi ) = {\psi}_1(s) \xi + {\psi}_2(s)(1 - \xi ).
\end{equation}

 In the transformed system we can use
the knowledge about the two-dimensional system where the variable
$\xi$ is constant. The systems, where $\xi=0$ and $\xi=1$ and where the 
three dimensional system has behaviour related to the two dimensional 
systems in the coordinate planes $x_1=0$ and $x_2=0$ correspondingly,
turn out to be important as comparison systems.

\subsection{Comparison systems}

This subsection contains the main results obtained by using the
comparison systems mentioned in previous subsection.
From these it is easy to construct a positive invariant set.

We will need some notations:
 $$H = H(x,s) = h(s)-{\psi}_1(s)x_1 - {\psi}_2(s)x_2,\qquad
 H_i=H_i(x,s)=h(s)-m{\psi}_i(s),$$

 \begin{equation}
 \omega(s) = {\psi}_1(s){\phi}_2(s) - {\psi}_2(s){\phi}_1(s) ,\qquad
 l=l(x,s)=h(s)({\phi}_2(s)-{\phi}_1(s)) - \omega(s)m .
 \end{equation}

 We define the following surfaces:
 \begin{equation}
 {\cal L} = \{ (x,s) \in R_3^+ \mid \; l(x,s)=0 \}\, ,  \qquad
 {\cal S}_i = \{ (x,s) \in R_3^+ \mid \; H_i(x,s)=0 \}.
 \end{equation}
 
 In the future we will use these notations both for the functions
    $\; l(m,s), H_i(m,s)\;$   and for the corresponding
    curves in the two-dimensional space
 $\,m,s. $
 \vskip1mm


  We consider two comparison systems 
  obtained for $\xi=1$ and $\xi=0$ and the corresponding vector
  fields
 ${\cal F}_i$:
 \begin{equation}\label{asr}
 \dot{x_1}={\phi}_i(s) x_1, \quad  \dot{x_2} = {\phi}_i(s)x_2,  \quad\quad
 \dot{s} = H_i ( x_1, x_2, s ).
 \end{equation}
 
  We notice that the planes, where $\xi$ 
 (and thus the proportion between the variables $x_1$ and $x_2$) is constant is invariant under these systems.
 
 These systems have integral surfaces generated by solutions in one plane.
 An integral surface is determined by an initial condition
 $m = m_0,\; s=s_0$. It contains this line and is usually a "spiral" surface, 
 but  a tube if the intial condition is on the limit cycle and the line
 if the initial condition is the equilibrium. This surface can be considered as having an outer surface with normal vector pointing away from the line of 
 equilibrium points. 
  We denote representants for the integral surfaces by 
${\cal{M}}_1$ and ${\cal{M}}_2$  correspondingly.

We  define the fields of  {\em outer normal vectors}
 $\overline{\bf n}_1$   and   $\overline{\bf n}_2$:
 \begin{equation}
 \overline{\bf n}_i =
 \overline{\bf n}_i\,( x_1, x_2, s ) =
 \,\left(\,-H_i(x_1, x_2, s),\,-H_i(x_1, x_2, s),\,  m\,\phi_i (s)\,\right).
 \end{equation}

 The vectors are orthogonal to solution curves of the comparison 
 systems and to the line $m = m_0,\; s=s_0.$ 
 Thus the vector fields
  $\overline{\bf n}_i$ degenerate at the equilibrium points 
 ${\cal F}_i $ and they are continuous and ortogonal 
 to integral surfaces of systems  $(\ref{asr}_i), $
 passing through line $m = m_0,\; s=s_0.$ 
 
 \vskip1mm

In order to find out in which direction the solutions of the
original three dimensional system with vector field ${\cal F}$ crosses
the integral surfaces ${\cal{M}}_1$ and ${\cal{M}}_2$,
we consider the scalar product:
 \begin{eqnarray}\label{11}
 <\,\overline{\bf n}_1,\,{\cal F}\,> \; = &\; x_2\,[ h(\phi_1 -\phi_2)
 + m\,(\psi_1 \phi_2 - \phi_1 \psi_2)\,]\,= -l(m,s)\,x_2, \\ \label{12}
 <\,\overline{\bf n}_2,\,{\cal F}\,> \; = &\; x_1\,[ h(\phi_2 -\phi_1)
 + m\,(\psi_2 \phi_1 - \phi_2 \psi_1)\,]\,= l(m,s)\,x_1.\\
 \nonumber
 \end{eqnarray}

We note that these products are zero on the surface $\cal L$.

For this surface to have a nice form we introduce
 a new assumption\\

 $B_1 :\quad$ $\omega (s)\neq 0$ for $s>0.$\\

From $A_3,\, A_5,\, A_6$ and $A_7$ follows that 
$\omega (s)>\phi_1(s) (\psi_2(s)-\psi_1(s))>0$ for $s>\lambda_1$
 and thus from $A_1$ and the assumption follows that  $\omega (s)>0$
for $s>0.$  

From this assumption also follows that, in region
 $s>0$, the surface
$\cal L$ can be defined explicitly:
\begin{eqnarray}
m=m^*(s)= h(s)\,\frac{\phi_2(s)-\phi_1(s)}{\omega(s)}.
 \end{eqnarray}

The surface is defined only for
$ \tau\leq s\leq 1$.

From sign analysis of expressions
 \ref{11} and \ref{12}, we get the following result:
 
\smallskip

{\bf Lemma 1.}  {\it 
The trajectories of system
 ($\ref{a1}$) intersect, in the region $m>m^*(s)$, transversally integral surfaces of
 ${\cal M}_1$ from its inner side to the outer and the integral surfaces 
  ${\cal M}_2$ from the outer side  to the inner side.

The intersection occurs in opposite order in the region $m<m^*(s)$.

The trajectories of the systems 
 (\ref{a1}) and (\ref{asr}$_i$)
are tangent to each other on the surface $m=m^*(s)$,
in the projection to the $(m,s)$-space. }

{\it Proof}. The type of intersection with integral surfaces follows from that the fact, that the first scalar product is positive and the second
negative in $m>m^*(s)$. The scalar products are zero for $m=m^*(s)$, from which 
follows the tangency. Thus the lemma is proved.

\smallskip

The systems
 (\ref{asr}$_i$) are projected to the  $(m,s)$-plane in a natural way:

 \begin{eqnarray}\label{14}
\dot{m}=  \phi_1(s) m,\quad & \dot{s} = h(s)-\psi_1(s) m; \\
\label{15}
\dot{m}=   \phi_2(s) m,\quad & \dot{s} = h(s)-\psi_2(s) m. \\
 \nonumber
 \end{eqnarray}

This defintion will give us one of the main technical
concepts of this work.

\smallskip

{\bf Definition 1}.
The projection of the surface
 $x+y=m^*(s)$ into the  $(m,s)$-plane forms a  
curve $m=m^*(s)$, we will call {\sf tangency curve.}

\smallskip

We also need the following definition.

\smallskip

{\bf Defintion 2}.
The projections $S_1$ and $S_2$ of the isoclines of the comparison systems,
defined by the equalities
$m=h(s)/\psi_1(s)$ \ and  \  $m=h(s)/\psi_2(s)$ 
we will call  {\sf first and second main
isoclines}
correspondingly.

\smallskip

The following lemma gives us information about the position of 
the tangency curve and the main isoclines.

\smallskip

{\bf Lemma 2}.
 {\it The main isoclines do not intersect in the positive quadrant of the
 $(m,s)$-space. Thereby the 
 first main isocline $S_1$ will be outside 
 the second $S_2$ relative to $s$-axis.

The tangency curve intersects the main isoclines at equilibrium 
points of systems
 (\ref{14}-\ref{15}) 
and is outside the first main isocline, if 
 $s>\lambda_1,$  between the main isoclines, if 
 $\lambda_2<s<\lambda_1,$ and inside the second, if
$s<\lambda_2.$  }

\smallskip

{\it Proof}. The position of the main isoclines follows from $A_5$. The inequality $m^*(s) > \frac{\displaystyle h(s)}{\displaystyle\psi_i(s)}$ is equivalent
to 
$\phi_2(s) (\psi_i(s)-\psi_1(s))>\phi_1(s) (\psi_i(s)-\psi_2(s)) $ because $\omega (s)>0$.
But then from $A_5$ follows equivalence with $\phi_i(s)>0$,
 which follows from 
 $s>\lambda_i$ using  $A_3$. Analogously 
 $m^*(s)< \frac{\displaystyle h(s)}{\displaystyle\psi_i(s)}$
follows from $s<\lambda_i$. These inequalities give the geometric conclusion of the lemma.

 \subsection{The standard system}
 We examine more carefully the expressions for the tangency curve in the case where we use  (\ref{a2})
 and (\ref{a3}) in the right hand sides.
 
 We calculate:

 \begin{equation}\label{16}\phi_2(s)-\phi_1(s) =\frac{\kappa_0
(s-\tau)}{(s+a_1)(s+a_2)},\quad \mbox{where} \quad
 \kappa_0 =a_1+\lambda_1 -(a_2+\lambda_2), \
 \tau=\frac{(\gamma-1)\lambda_1a_2}{\kappa_0}.
 \end{equation}

Moreover,

 \begin{equation}\label{17}
\omega(s)=-\frac{s-\lambda_1}{s+a_1}\cdot\frac{s}{s+a_2}+
\frac{s-\lambda_2}{s+a_2}\cdot\frac{s}{s+a_1}
=\frac{s(\lambda_1-\lambda_2)}{(s+a_1)(s+a_2)}.
 \end{equation}

\noindent
and the formula for the tangency curve becomes a quadratic expression:

 \begin{equation}\label{18}
m^*(s)=\frac{h(s)\kappa_0(s-\tau)}{s(\lambda_1-\lambda_2)}
 = (1-s)(s-\tau)\kappa, \quad \mbox{where} \quad \kappa =
 \frac{\kappa_0}{\lambda_1-\lambda_2}.
 \end{equation}

We notice that parameters giving the same tangency curve satisfy the equation

\begin{equation}
a_i=\lambda_i (\kappa -1) -\tau\kappa
\label{equ}
\end{equation}

 \subsection{Construction of the boundaries 
 of the Poincar\'{e} annulus}
 
   Using the lemmas we construct an annular region where trajectories
  cross  the boundaries strictly in the inside direction. We use the region
  to define a Poincar\'{e} map to discretize the flow inside the
  annulur region. This region we will call 
 {\sf Poincar\'{e} annulus.}

  We will build it in the
  $(m,s)$-plane. Its preimage in the  
   $(m,s,\xi)$-space is a positively invariant
   region possible containing a non-wandering set
   with complicated structure.
   Also this 2-dimensional region we call Poincar\'{e} annulus.

We define  two continuous systems by gluing together  
 two parts:

1) System ${\cal P}^i$ is defined by the vector field, which 
coincides with the first comparison system (\ref{14}) outside 
the tangency curve
 $\cal L$ (that's in region $m>m^*(s))$ and with the second
 system  (\ref{15}) inside. We will call this system 
the {\sf inner} system.

2) System ${\cal P}^e$ is defined opposite:
the vector field coincides with the second comparison system (\ref{15})
outside the tangency curve $\cal L$ and with the first
system (\ref{14}) inside. 
 We will call this system 
the {\sf outer} system.

The orbits of the systems are correctly defined 
and continuous
 everywhere except
on the segment  $[O_2,O_1]$ of the tangency curve.
This follows from the fact that the comparison vector fields
are collinear and non-degenerate on the curve 
 $\cal L,$ and intersect this curve transversally.

The segment  $[O_2,O_1]$ we denote by 
 $\cal A$. 

For the outer system the segment $\cal A$ is repelling
(orbits with initial conditions not on   $\cal A$ never
come on  $\cal A$), and for the inner system this segment is
attracting. The outer  system has one stable limit cycle and the inner system
has two limit cycles of which the outer is stable and the inner unstable.

\paragraph{Construction of outer boundary.} 
We look at the separatix of the unstable manifold of the equilibrium point
 $O_1$ of the outer system ${\cal P}^e$ and we denote its first
 intersection with the tangency curve $\cal L$, after passing line $s=\lambda_1$,
 by  $L_1$. It is clear that the intersection is below the
 segment $\cal A$. The next intersection of the orbit
 with the tangency curve we denote by   $L_2.$
 Thereby we have found a positively invariant set for the system. The boundary of
 it is given by the part of the separatrix and of $\cal L$ from $O_1$ to $L_2$ and   we will  take it as one boundary of the Poincar\'{e} annulus.

\paragraph{Construction of inner boundary.}
  The outer boundary of the Poincar\'{e} annulus was constructed 
  without problems. But the construction of the inner boundary
  is more problematic.
  We explain why.  
  We choose an arbitrary point  $M$ on the curve $\cal L$ as the initial point
  for on orbit of the inner system ${\cal P}^i$. 
  We will suppose that the initial point coincides with the
  upper point on $\cal L$ of an orbit of ($\ref{15}$) tending to the limit cycle 
   from the
  inner side.
   There are two possibilities.  
    Either the orbit of the point $M$ does not intersect with
     the segment  $\cal A$ or it intersects the segment.
     
\smallskip

     In the case of non-intersection the orbit tends to a limit cycle
     of system ${\cal P}^i$, which we will take
     as inner boundary of the Poincar\'{e} annulus.

  In this case, the Poincar\'{e} map for the full system
 (\ref{a1}) is defined in the following way:

In the $(m,s,\xi)$-space (or correspondingly in the 
$(x,y,s)$-space) we consider the plane
 $S_{\varepsilon}: \, s=\varepsilon,$ where 
$\varepsilon$ 
is chosen such that the intersection of the plane
with the Poincar\'{e} annulus consists of two rectangles
in the $(m,s,\xi)$-space. On $S_{\varepsilon}$ the  Poincar\'{e} map is defined as usually as the next
 intersection with $S_{\varepsilon}$.
 
In this case we will say that the Poincar\'{e} map is correctly defined.

\smallskip

In the case, when the closure of the orbit has common points 
with the special segment, we will speak about non-correct boundary
for the Poincar\'{e} annulus.

If the Poincar\'{e} map is correctly defined for some parameter values in the
standard system given by functions (\ref{a2}) and (\ref{a3}) and we increase $\lambda_2$ and
let $a_2$ satisfy the equality $(\ref{equ})$, it remains correctly defined.
This holds, because if we take the limit cycle for the inner boundary with 
original parameters and consider the trajectory containing the outer part
of this cycle (outside tangency curve, where trajectories for new system coincide
with original) for the new $\lambda_2$ and $a_2$ then from Lemma 1 follows  that the trajectory for the new system inside tangency curve will be outside the limit
cycle.

Some concrete numerical results for the boundary are given in next section.

\begin{figure}[h!]
\begin{center}
\caption{Illustration of tangency curve, main isoclines and construction
of boundary for outer system}
\includegraphics[height=7.5cm,width=15cm]{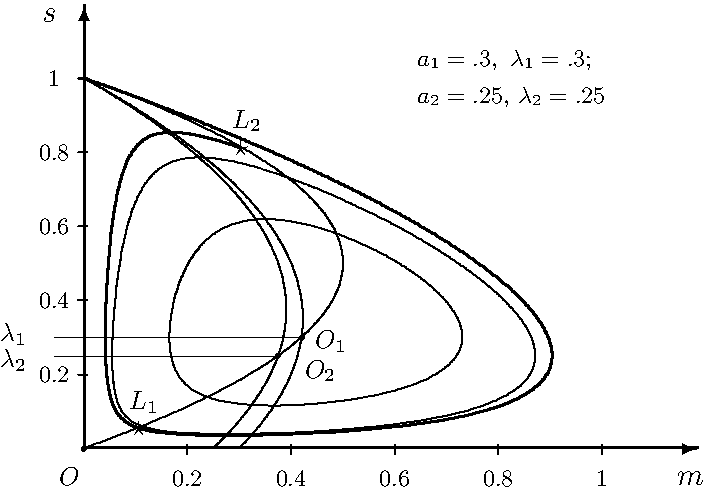}\\
\label{poi2}
\end{center}
\end{figure}

\section{Numerical results}

All results here are given for the standard system given by functions (\ref{a2})  and
(\ref{a3}).

\begin{figure}[h!]
\begin{center}
\caption{Bifurcation diagram when $a_1=\lambda_1=0.1$ and $\lambda_2=\nu \, 0.01$ and
$a_2$ is chosen so that we stay on tangency curve determined by 
$a_2=0.0075,\, \lambda_2=0.01$. Horisontal axis represents the parameter $\nu$ and
vertical represents the variable $\xi$. The Poincar\'{e} map is defined on the
outer rectangle of the Poincar\'{e} annulus on level $s=0.1$}
\includegraphics[height=7.5cm,width=15cm]{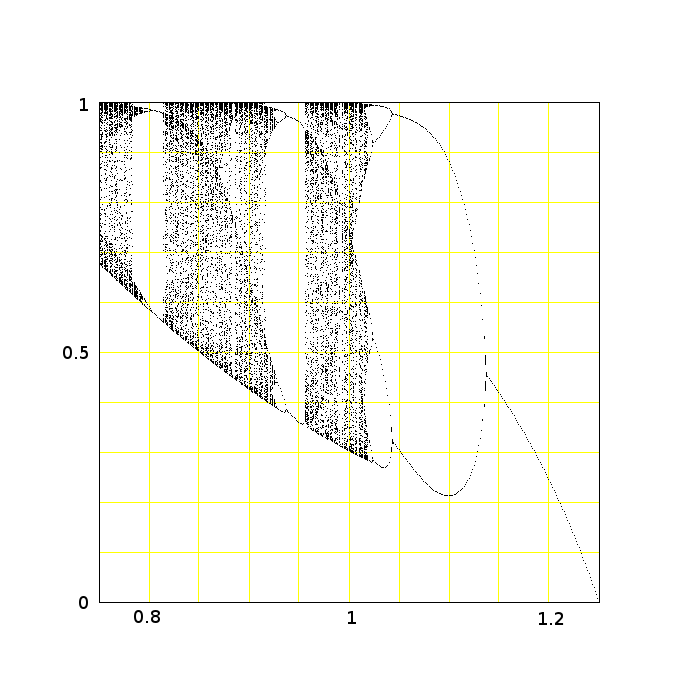}\\
\label{gunfig}
\end{center}
\end{figure}

We choose $a_1=\lambda_1=0.1, \, a_2=0.0075,\,\lambda_2=0.01$ and calculate 
$\tau\approx 0.00137, \, \kappa \approx 2.028$. 
In this case the stable limit cycle of the inner system is intersecting
the tangency curve at about $s=0.85$ and the unstable cycle at $0.384$.
If we increase $\lambda_2$ to 0.02 and calculate $a_2$ from formula
(\ref{equ}), the stable limit cycle of the inner system is intersecting
the tangency curve at about $s=0.905$ and the unstable cycle at $s=0.2703$. 
If we increase further $\lambda_2$ to 0.05 and calculate $a_2$ from formula
(\ref{equ}), the stable limit cycle of the inner system is intersecting
the tangency curve at about $s=0.95$ and the unstable cycle at $s=0.1595$. 
A bifurcation diagram showing how attractor for Poincar\'{e} map changes with
$\lambda_2$ is given in Figure \ref{gunfig}.

In the case when the Poincar\'{e} map is correctly defined, very often there is a strong
contraction in the $m$-direction and in \cite{osipov3} is shown by numerical experiments and theoretical estimating arguments that the one dimensional model map
given by

$$ f(v)= \beta + v -\frac{k_1+k_2 e^ v}{1+v} u$$ 

\noindent where  $\beta ,\, u$ and $k_i$ are constants and $v=\ln (x_2/x_1)$
gives a good approximation.

In the case when the Poincar\'{e} map is not defined correctly we can find
interesting "spiral-like" chaos for parameter values
$a_1=0.5,\, \lambda_2=0.2,\, a_2=0.002$  and $0.33\leq \lambda_1\leq 0.42$ 
and $a_1=0.2,\, \lambda_2=0.2,\, a_2=0.01$  and $0.4\leq \lambda_1\leq 0.5$.
We observe that in this case the equilibrium in $x_1 s$-plane is stable
inside the coordinate plane and there is no periodic solution.
Examples of more  complicated chaos are also easy to find when there 
are limit cycles in both coordinate planes, for example, if
$a_1=\lambda_1 =0.25,\, \lambda_2=0.163,\, a_2=0.06$.

\section{Conclusion}
We have introduced the concept of tangency curve as an important tool for constructing a positively invariant set, where we can define a 
Poincar\'{e} map for examining the structure of attractors of system (\ref{a1}).
For the special case, where functions are given by (\ref{a2}) and (\ref{a3}), 
we have given a concrete numerical example of a bifurcation diagram for this Poincar\'{e} map.
Except general conditions on the system, we have introduced a geometric condition for the existence of our invariant set. Our conditions are sufficient for the existence of the invariant set, but not necessary. For example, the bifurcation path in \cite{osipov4} does not satisfy the geometric condition for all parameters, but it is possible to find a good invariant set, where to define a Poincar\'{e} map.
In general, in case of functions (\ref{a2}) and (\ref{a3}), there is a lower limit for parameter $\lambda_2$ when our conditions are not more satisfied. We also give numerical examples of more complicated "spiral-like" case where no similar invariant set of the type we constructed is possible to find.

\bibliographystyle{plain}
\bibliography{osipPoinc}

 \end{document}